\documentclass{article}
\usepackage{amsthm}
\usepackage{amssymb}
\usepackage[applemac]{inputenc}

\newtheorem{theorem}{Theorem}
\newtheorem{corollary}{Corollary}
\newtheorem{definition}{Definition}

\def\pts{\textup{\,:\,}}

\title{Independence of the metric in the fine $C^0$-topology of a function space}
\author{Francisco J. Gonzalez-Acu\~na}

\begin{document}
\maketitle

\noindent
The following is a translation of the paper ``{\em Independencia de la métrica en la $C^0$ topología fina de un
espacio de funciones}''. Vol.\ 3 (1963) pp.\ 29-35.\ Anales del Instituto de Matemáticas (Universidad Nacional Autónoma
México). I thank Francisco Marmolejo for writing this translation.

\maketitle

\section{Introduction}

Several topologies can be given to the space of continuous functions $F(X,Y)$ from a topological
space $X$ into a metrizable space $Y$. One of them is defined as follows: let $d$ be a bounded distance
consistent with the topology of $Y$. A distance $d^*$ can be defined on $F(X,Y)$ by:
\[
d^*(f,g)=\sup_{x\in X}d[f(x),g(x)].
\]
The topology of $F(X,Y)$ determined by $d^*$ is called the {\em $d^*$-topology}. This topology depends not
only on the topologies of $X$ and $Y$ but also on the metric $d$ that is chosen for $Y$.

In \S\ref{sec2} we define de fine $C^0$-topology on $F(X,Y)$. It is shown in [1,\S5] that this topology does not depend
on the chosen metric for $Y$ in the case where $X$ paracompact. We show here (Theorem \ref{teo1})
the same result without imposing any condition on $X$. 

In \S\ref{sec3} some properties of the fine $C^0$-topology are mentioned and it is compared with other topologies.

\section{The fine $C^0$-topology}\label{sec2}

\begin{definition}
Let $X$ be a topological space, $(Y,d)$ metric, $f\pts X\to Y$ continuous, $\delta\pts X\to\mathbb{R}^+=\{$positive reals$\}$
continuous. $g\pts X\to Y$ continuous is a $\delta$-approximation of $f$ if $d[f(x),g(x)]<\delta(x)$ for all $x\in X$.
The {\em $\delta$-neighborhood} of $f$ is the set of all the $\delta$-approximations of $f$. If $F(X,Y)$ is the set of
continuous functions from $X$ to $Y$ and $f\in F(X,Y)$, we define a neighborhood of $f$ as a subset of $F(X,Y)$
that contains some $\delta$-neighborhood of $f$.  
\end{definition}

This defines a topology on $F(X,Y)$ given that the intersection of the 
$\delta$-neighborhood of $f$ and the $\eta$-neighborhood of $f$ is the $\min(\delta,\eta)$-neighborhood of $f$ and 
furthermore if $g$ is in the $\delta$-neighborhood of $f$ and $\eta(x)=\delta(x)-d[f(x),g(x)]$, then the $\eta$-neighborhood 
of $g$ is contained in the $\delta$-neighborhood of $f$. We will call this topology the {\em fine $C^0$-topology} of 
$F(X,Y)$, and we denote by $T_d(X,Y)$ the resulting topological space.

The next theorem establishes that the topology on $T_d(X,Y)$ depends only on the topologies of $X$ and $Y$ and not
on the metric on $Y$.

\begin{theorem}\label{teo1}
If $d_1$ and $d_2$ are two equivalent distances on $Y$, then $T_{d_1}(X,Y)$ and $T_{d_2}(X,Y)$ are the same 
topological space.
\end{theorem}

\proof Let $W_1$ be a $\delta_1$-neighborhood of $f\pts X\to Y$ on $T_{d_1}(X,Y)$. We will show that there exists
$W_2$, $\delta_2$-neighborhood of $f$ in $T_{d_2}(X,Y)$, contained in $W_1$. The proof that for every
$\delta_2$ neighborhood of $W_2$ of $f$ there is a $\delta_1$-neighborhood contained in $W_2$ being
analogous. 

We define $\delta_2\pts X\to \mathbb{R}^+$ as
\[
\delta_2=G'\circ(f\times \delta_1),
\]
where $G'$ is a continuous function from $Y\times \mathbb{R}^+$ to $\mathbb{R}^+$ that satisfies:
\[
G'(y,t)\leq \sup\{r|B_2(y,r)\subset B_1(y,t)\}=G(y,t) 
\]
($B_i(y,r)=\{y'\in Y|d_i(y,y')<r\}$ $i=1,2$). Let's show that there is such a $G'$.

The function $G$, that is positive since $d_1$ and $d_2$ are equivalent, satisfies:
\[
G(y',t)\geq G(y,t-\varepsilon)-\varepsilon \textup{ \ if \ } d_i(y,y')<\varepsilon, \ \ i=1,2.
\]
Indeed, if $w\in B_2(y',G(y,t-\varepsilon)-\varepsilon)$ and $d_i(y,y')<\varepsilon$ $i=1,2$,\medskip

$d_2(w,y')+\varepsilon < G(y,t-\varepsilon)$\medskip

$d_2(w,y)<G(y,t-\varepsilon)$\medskip

$B_2(y,s)\subset B_1(y,t-\varepsilon) \textup{ \ with \ } d_2(w,y)<s<G(y,t-\varepsilon)$\medskip

$d_1(w,y)<t-\varepsilon$\medskip

{\ }\hspace{2.6cm} $w\in B_1(y',t)$\medskip

and\medskip

$B_2(y',G(y,t-\varepsilon)-\varepsilon)\subset B_1(y',t)$\medskip

that is\medskip

$G(y',t)\geq G(y,t-\varepsilon)-\varepsilon$.\medskip

Since $G$ is not decreasing in the second variable, for $(y,t)\in Y\times\mathbb{R}^+$ there is $\varepsilon$
such that $G(y,t-\varepsilon)-\varepsilon>0$ and thus every point $(y,t)$ of $Y\times\mathbb{R}^+$ has a
neighborhood in which $G$ has a positive lower bound, namely the neighborhood 
$[B_1(y,\varepsilon)\cap B_2(y,\varepsilon)]\times (a,\infty)$ with $0<a<t$ and $G(y,a-\varepsilon)-\varepsilon>0$.

Hence, since the domain of $G$ is paracompact, (every metrizable is paracompact), there is a locally finite open cover
$\{V_\alpha\}$ of this domain such that $G$ has a positive lower bound $\varepsilon_\alpha$ on each $V_\alpha$.
Let $\{\phi_\alpha\}$ be a partition of unity subordinated to $\{V_\alpha\}$. Define 
$G'(y,t)=\sum_\alpha\varepsilon_\alpha\phi_\alpha(y,t)$, a continuous function with values in $\mathbb{R}^+$.

We have that
\[
G'(y,t)=\sum_{(y,t)\in V_\alpha}\varepsilon_\alpha\phi_\alpha(y,t) \leq
\max_{(y,t)\in V_\alpha}\{\varepsilon_\alpha\} \leq G(y,t).
\]
If $W_2$ is the $\delta_2$-neighborhood of $f$ in $T_{d_2}(X,Y)$, $W_2\subset W_1$ since, if $g\in W_2$,
\medskip

$d_2[g(x),f(x)]<\sup\{r|B_2(f(x),r)\subset B_1(f(x),\delta_1(x))\}$\medskip

$B_2(f(x),s)\subset B_1(f(x),\delta_1(x))$ with\medskip

$d_2[g(x),f(x)]<s<\sup\{r|B_2(f(x),r)\subset B_1(f(x),\delta_1(x))\}$\medskip

$d_1[f(x),g(x)]<\delta_1(x)$\medskip

and\medskip

$g\in W_1$\medskip

This completes the proof. 
\medskip

We can thus omit the subindex $d$ in $T_d(X,Y)$.

\section{Properties and comparison with other \\ topologies}\label{sec3}

From now on $Y$ will always denote a metrizable space. 

The space $T(X,Y)$ is always Tychonoff ($T_1$ and completely regular) since if $f\in T(X,Y)$ and
$\delta$ is any positive continuous function, the function from $T(X,Y)$ to the reals defined by
\[
g\mapsto \min\left[\sup_{x\in X}\left\{\frac{d[f(x),g(x)]}{\delta(x)}\right\},1\right],
\]
where $d$ is a distance in $Y$ consistent with the topology of $Y$, is continuous, has value 1 outside of the
$\delta$-neighborhood of $f$, and has value $0$ if and only if $g=f$.

A topology in $F(X,Y)$ is called {\em admissible} if the function $(f,x)\mapsto f(x)$ defined on $F(X,Y)\times X$ turns
out to be continuous when $F(X,Y)$ is endowed with that topology. The fine $C^0$-topology and the $d^*$
topology are admissible. The compact-open topology is coarser that the $d^*$-topology, which in turn is coarser
than the fine $C^0$-topology [2].

In case $X$ is compact these topologies are identical (see for example [3,Chap.\ 7, Theorem 11]). We will see that
if $X$ is $T_4$, not countably compact and $Y$ has a subspace homeomorphic to the reals, then the 
three topologies are all different. 

\begin{theorem}\label{teo2}
If $X$ is countably compact (every open countable covering has a finite subcovering), the fine $C^0$-topology of
$F(X,Y)$ coincides with the $d^*$-topology. If $X$ is $T_1$ and normal, $Y$ has a subspace homeomorphic to
$\mathbb{R}$ and $T(X,Y)$ satisfies the first countability axiom, then $X$ is countably compact.
\end{theorem}

\proof To show the first statement it suffices to show that every positive continuous function $\delta$ on a
countably compact space $X$ has a positive lower bound.

If $X$ is countably compact, $\delta(X)$ is as well. A countably compact in $\mathbb{R}$ is compact, so
$\delta(X)$ has a positive lower bound. 

Assume now that $X$ is $T_1$, normal and not countably compact, and
that $Y$ contains $\mathbb{R}$ as a subspace. We will show that $T(X,Y)$ does not satisfy the first countability
axiom. 

$T(X,\mathbb{R})$ is a subspace of $T(X,Y)$ so it will suffice to show that $T(X,\mathbb{R})$ does not satisfy
the first countability axiom.

Since $X$ is not countably compact there is a sequence of distinct points $(x_i)$ without a cluster point in $X$. 
$\{x_i\}$ is closed in $X$ and discrete since $X$ is a $T_1$ space.

Let $f\in T(X,\mathbb{R})$ be the function identically $0$, and let $\{\delta_i\}_{i=1,2,\dots}$ be any countable family
of positive continuous functions on $X$.  By the normality of $X$ we can define $\delta\pts X\to\mathbb{R}^+$
continuous such that $\delta(x_i)=\frac 12\delta_i(x_i)$, $i=1,2,\dots$.

The $\delta$-neighborhood of $f$ does not contain any $\delta_i$-neighborhood of $f$. Indeed, 
since $X$ is normal, for every $i$ there is $g\in T(X,\mathbb{R})$ such that $g(x_i)=\delta_i(x)$,
$g(X)=[0,\delta(x_i)]$ and $g(X-V)=\{0\}$ where $V$ is a neighborhood of $x_i$ in which $\delta_i$ is bigger
than $\frac 12\delta_i(x_i)$. $g$ is then in the $\delta_i$-neighborhood of $f$ but not in the $\delta$-neighborhood
of $f$ (taking the usual distance in $\mathbb{R}$).

Therefore the $\delta_i$-neighborhoods of $f$ do not form a fundamental system of neighborhoods of $f$.
This completes the proof of the theorem.\medskip

Condition $T_1$ can not be omitted from the statement of the theorem since, if $X$ is the set of natural numbers
with the topology in which the open sets are the sets of the form $\{1,2,\dots,n\}$, the empty set and the whole $X$, then
$X$ is normal, not countably compact and $T(X,Y)$ satisfies the first countability axiom.

The proposition is also not valid if normal is replaced by completely regular. The following is a counterexample:
$X=\Omega'\times \omega'-\{(\Omega,\omega)\}$ where $\Omega'$ and $\omega'$ are the set of ordinals not
bigger that the first  uncountable ordinal $\Omega$ and the set of ordinals not bigger than the first infinite ordinal
$\omega$ respectively, both spaces with the order topology. $X$ is $T_1$, completely regular, not countably
compact and $T(X,Y)$ satisfies the first countability axiom.

\begin{corollary}
If $X$ is $T_1$, normal, not countably compact and $Y$ has a subspace homeomorphic to $\mathbb{R}$, then the
compact-open topology, the $d^*$-topology and the fine $C^0$-topology on $F(X,Y)$ are all different. See 
\textup{[2,\S3]} 
\end{corollary}

\section*{References}

\begin{tabular}{lp{8cm}}
{[1]}  Whitehead, J.H.C.: & Manifolds with transverse fields in euclidean space, Ann.\ of Math.\ vol.\ 73 (1961), 
pp 154-212.\medskip\\
{[2]}  Jackson, J.R.: & Comparison of topologies on function spaces. Proc.\ Amer.\ Math.\ Soc.\ vol.\ 3 (1952), 
pp 156-158.\medskip\\
{[3]} Kelley, J. L.: & General topology (New York, 1955).
\end{tabular}

\end{document}